# CONVERGENCE OF ADAPTIVE MIXTURES OF IMPORTANCE SAMPLING SCHEMES[1]

By R. Douc, A. Guillin, J.-M. Marin and C. P. Robert

*École Polytechnique, École Centrale Marseille et LATP, CNRS, INRIA Futurs, Projet Select, Université d'Orsay and Université Paris Dauphine and CREST-INSEE*

In the design of efficient simulation algorithms, one is often beset with a poor choice of proposal distributions. Although the performance of a given simulation kernel can clarify a posteriori how adequate this kernel is for the problem at hand, a permanent on-line modification of kernels causes concerns about the validity of the resulting algorithm. While the issue is most often intractable for MCMC algorithms, the equivalent version for importance sampling algorithms can be validated quite precisely. We derive sufficient convergence conditions for adaptive mixtures of population Monte Carlo algorithms and show that Rao–Blackwellized versions asymptotically achieve an optimum in terms of a Kullback divergence criterion, while more rudimentary versions do not benefit from repeated updating.

## 1. Introduction.

1.1. *Monte Carlo calibration.* In the simulation settings found in optimization and (Bayesian) integration, it is well documented [20] that the choice of the instrumental distributions is paramount for the efficiency of the resulting algorithms. Indeed, in an importance sampling algorithm with importance function $g(x)$, we are relying on a distribution $g$ that is customarily difficult to calibrate outside a limited range of well-known cases. For instance, a standard result is that the optimal importance density for approximating an integral

$$\mathfrak{I} = \int f(x)\pi(x)\,dx$$

Received January 2005; revised December 2005.

[1]Supported in part by an ACI "Nouvelles Interfaces des Mathématiques" grant from the Ministère de la Recherche.

is $g^{\star}(x) \propto |f(x)| \pi(x)$ ([20], Theorem 3.12), but this formal result is not very informative about the practical choice of $g$, while a poor choice of $g$ may result in an infinite variance estimator. MCMC algorithms somehow attenuate this difficulty by using local proposals like random walk kernels, but two drawbacks of these proposals are that they may take a long time to converge [18] and their efficiency ultimately depends on the scale of the local exploration.

The goals of Monte Carlo experiments are multifaceted and therefore the efficiency of an algorithm can be evaluated from many different perspectives. In particular, in Bayesian statistics the Monte Carlo sample can be used to approximate a variety of posterior quantities. Nonetheless, if we try to assess the generic efficiency of an algorithm and thus develop a portmanteau device, a natural approach is to use a measure of agreement between the target and the proposal distribution, similar to the intrinsic loss function proposed in [19] for an invariant estimation of parameters. A robust measure of similarity ubiquitous in statistical approximation theory [7] is the Kullback divergence

$$\mathfrak{E}(\pi, \tilde{\pi}) = \int \log \frac{d\pi(x)}{d\tilde{\pi}(x)} \pi(dx),$$

and this paper aims to minimize $\mathfrak{E}(\pi, \tilde{\pi})$ within a class of proposals $\tilde{\pi}$.

1.2. *Adaptivity in Monte Carlo settings.* Given the complexity of the original optimization or integration problem (which does itself require Monte Carlo approximations), it is rarely the case that the optimization of the proposal distribution against an efficiency measure can be achieved in closed form. Even the computation of the efficiency measure for a given proposal is impossible in the majority of cases. For this reason, a number of adaptive schemes have appeared in the recent literature ([20], Section 7.6.3) in order to design better proposals against a given measure of efficiency without resorting to a standard optimization algorithm. For instance, in the MCMC community, sequential changes in the variance of Markov kernels have been proposed in [13, 14], while adaptive changes taking advantage of regeneration properties of the kernels have been constructed by Gilks, Roberts and Sahu [12] and Sahu and Zhigljavsky [23, 24]. From a more general perspective, Andrieu and Robert [2] have developed a two-level stochastic optimization scheme to update parameters of a proposal towards a given integrated efficiency criterion such as the acceptance rate (or its difference with a value known to be optimal—see Roberts, Gelman and Gilks [21]). However, as reflected in this general technical report of Andrieu and Robert [2], the complexity of devising valid adaptive MCMC schemes is a genuine drawback in their extension, given that the constraints on the inhomogeneous Markov chain that results from this adaptive construction are either difficult to satisfy or result in a fixed proposal after a certain number of iterations.

Cappé, Guillin, Marin and Robert [3] (see also [20], Chapter 14) developed a methodology called Population Monte Carlo (PMC) [16] motivated by the observation that the importance sampling perspective is much more amenable to adaptivity than MCMC, due to its unbiased nature: using sampling importance resampling, any given sample from an importance distribution $g$ can be transformed into a sample of points marginally distributed from the target distribution $\pi$ and Cappé et al. [3] (see also [8]) showed that this property is also preserved by repeated and adaptive sampling. (In this setting, "adaptive" is to be understood as a modification of the importance distribution based on the results of previous iterations.) The asymptotics of adaptive importance sampling are therefore much more manageable than those of adaptive MCMC algorithms, at least at a primary level, if only because the algorithm can be stopped at any time. Indeed, since every iteration is a valid importance sampling scheme, the algorithm does not require a burn-in time. Borrowing from the sequential sampling literature [11], the methodology of Cappé et al. [3] thus aimed at producing an adaptive importance sampling scheme via a learning mechanism on a population of points, themselves marginally distributed from the target distribution. However, as shown by the current paper, the original implementation of Cappé et al. [3] may suffer from an asymptotic lack of adaptivity that can be overcome by Rao–Blackwellization.

1.3. *Plan and objectives.* This paper focuses on a specific family of importance functions that are represented as mixtures of an arbitrary number of fixed proposals. These proposals can be educated guesses of the target distribution, random walk proposals for local exploration of the target, nonparametric kernel approximations to the target or any combination of these. Using these fixed proposals as a basis, we then devise an updating mechanism for the weights of the mixture and prove that this mechanism converges to the optimal mixture, the optimality being defined here in terms of Kullback divergence. From a probabilistic point of view, the techniques used in this paper are related to techniques and results found in [4, 6, 17]. In particular, the triangular array technique that is central to the CLT proofs below can be found in [4, 10].

The paper is organized as follows. We first present the algorithmic and mathematical details in Section 2 and establish a generic central limit theorem. We evaluate the convergence properties of the basic version of PMC in Section 3, exhibiting its limitations, and show in Section 5 that its Rao–Blackwellized version overcomes these limitations and achieves optimality for the Kullback criterion developed in Section 4. Section 6 illustrates the practical convergence of the method on a few benchmark examples.

**2. Population Monte Carlo.** The Population Monte Carlo (PMC) algorithm introduced in [3] is a form of iterated sampling importance resampling (SIR). The appeal of using a repeated form of SIR is that previous samples are informative about the connections between the proposal (importance) and the target distributions. We stress from the outset that this scheme has very few connections with MCMC algorithms since (a) PMC is not Markovian, being based on the whole sequence of simulations and (b) PMC can be stopped at any time, being validated by the basic importance sampling identity ([20], equation (3.9)) rather than a probabilistic convergence result like the ergodic theorem. These features motivate the use of the method in setups where off-the-shelf MCMC algorithms cannot be of use. We first recall basic Monte Carlo principles, mostly to define notation and to make precise our goals.

2.1. *The Monte Carlo framework.* On a measurable space $(\Omega, \mathcal{A})$, we are given a *target*, that is, a probability distribution $\pi$ on $(\Omega, \mathcal{A})$. We assume that $\pi$ is dominated by a reference measure $\mu$, $\pi \ll \mu$, and also denote by $\pi(dx) = \pi(x)\mu(dx)$ its density. In most settings, including Bayesian statistics, the density $\pi$ is known up to a normalizing constant, $\pi(x) \propto \tilde{\pi}(x)$. The purpose of running a simulation experiment with the target $\pi$ is to approximate quantities related to $\pi$, such as intractable integrals

$$\pi(f) = \int f(x)\pi(dx),$$

but we do not focus here on a specific quantity $\pi(f)$. In this setting, a standard stochastic approximation method is the Monte Carlo method, based on an i.i.d. sample $x_1, \ldots, x_N$ simulated from $\pi$, that approximates $\pi(f)$ by

$$\hat{\pi}_N^{\mathrm{MC}}(f) = N^{-1} \sum_{i=1}^{N} f(x_i),$$

which almost surely converges to $\pi(f)$ (as $N$ goes to infinity) by the law of large numbers (LLN). The central limit theorem (CLT) implies, in addition, that if $\pi(f^2) = \int f^2(x)\pi(dx) < \infty$, then

$$\sqrt{N}\{\hat{\pi}_N^{\mathrm{MC}}(f) - \pi(f)\} \overset{\mathscr{L}}{\rightsquigarrow} \mathcal{N}(0, \mathbb{V}_\pi(f)),$$

where $\mathbb{V}_\pi(f) = \pi([f - \pi(f)]^2)$. Obviously, this approach requires a direct i.i.d. simulation from $\pi$ (or $\tilde{\pi}$), which often is impossible. An alternative ([20], Chapter 3) is to use importance sampling, that is, to choose a probability distribution $\nu \ll \mu$ on $(\Omega, \mathcal{A})$ called the *proposal* or *importance* distribution, the density of which is also denoted by $\nu$, and to estimate $\pi(f)$ by

$$\hat{\pi}_{\nu,N}^{\mathrm{IS}}(f) = N^{-1} \sum_{i=1}^{N} f(x_i)\left(\frac{\pi}{\nu}\right)(x_i).$$

If $\pi$ is also dominated by $\nu$, $\pi \ll \nu$, then $\hat{\pi}_{\nu,N}^{\mathrm{IS}}(f)$ almost surely converges to $\pi(f)$ and if $\nu(f^2(\pi/\nu)^2) < \infty$, then the CLT also applies, that is,

$$\sqrt{N}\{\hat{\pi}_{\nu,N}^{\mathrm{IS}}(f) - \pi(f)\} \overset{\mathscr{L}}{\rightsquigarrow} \mathcal{N}\left(0, \mathbb{V}_\nu\left(f\frac{\pi}{\nu}\right)\right).$$

As the normalizing constant of the target distribution $\pi$ is unknown, it is not possible to directly use the IS estimator $\hat{\pi}_{\nu,N}^{\mathrm{IS}}(f)$. A convenient substitute is the self-normalized IS estimator,

$$\hat{\pi}_{\nu,N}^{\mathrm{SNIS}}(f) = \left(\sum_{i=1}^N (\pi/\nu)(x_i)\right)^{-1} \sum_{i=1}^N f(x_i)(\pi/\nu)(x_i),$$

which also converges almost surely to $\pi(f)$. If $\nu((1+f^2)(\pi/\nu)^2) < \infty$, then the CLT applies:

$$\sqrt{N}\{\hat{\pi}_{\nu,N}^{\mathrm{SNIS}}(f) - \pi(f)\} \overset{\mathscr{L}}{\rightsquigarrow} \mathcal{N}\left(0, \mathbb{V}_\nu\left\{[f - \pi(f)]\frac{\pi}{\nu}\right\}\right).$$

Obviously, the quality of the IS and SNIS approximations strongly depends on the choice of the proposal distribution $\nu$, which is delicate for complex distributions like those that occur in high-dimensional problems. (While multiplication of the number of proposals may offer some reprieve in this regard, we must stress at this point that our PMC methodology also suffers from the curse of dimensionality to which all importance sampling methods are subject in the sense that high-dimensional problems require a considerable increase in computational power.)

2.2. *Sampling importance resampling.* The sampling importance resampling (SIR) method of Rubin [22] is an extension of the IS method that achieves simulation from $\pi$ by adding resampling to simple reweighting. More precisely, the SIR algorithm operates in two stages. The first stage is identical to IS and consists in generating an i.i.d. sample $(x_1, \ldots, x_N)$ from $\nu$. The second stage builds a sample from $\pi$, $(\tilde{x}_1, \ldots, \tilde{x}_M)$, based on the instrumental sample $(x_1, \ldots, x_N)$, by resampling. While there are many resampling methods ([20], Section 14.3.5), the most standard (if least efficient) approach is multinomial sampling from $\{x_1, \ldots, x_N\}$ with probabilities proportional to the importance weights $[\frac{\pi}{\nu}(x_1), \ldots, \frac{\pi}{\nu}(x_N)]$, that is, the replacement of the weighted sample $(x_1, \ldots, x_N)$ by an unweighted sample $(\tilde{x}_1, \ldots, \tilde{x}_M)$, where $\tilde{x}_i = x_{J_i}$ $(1 \le i \le M)$ and where $(J_1, \ldots, J_M) \sim \mathcal{M}(M, \varrho_1, \ldots, \varrho_N)$, the multinomial distribution with probabilities

$$\mathbb{P}[J_\ell = i | x_1, \ldots, x_N] = \varrho_i \propto \frac{\pi}{\nu}(x_i), \qquad 1 \le i \le N, 1 \le \ell \le M.$$

The SIR estimator of $\pi(f)$ is then the standard average

$$\hat{\pi}_{\nu,N,M}^{\mathrm{SIR}}(f) = M^{-1} \sum_{i=1}^M f(\tilde{x}_i),$$

which also converges to $\pi(f)$ since each $\tilde{x}_i$ is marginally distributed from $\pi$. By construction, the variance of the SIR estimator is greater than the variance of the SNIS estimator. Indeed, the expectation of $\hat{\pi}_{\nu,N,M}^{\mathrm{SIR}}(f)$ conditional on the sample $(x_1, \ldots, x_N)$ is equal to $\hat{\pi}_{\nu,N}^{\mathrm{SNIS}}(f)$. Note, however, that an asymptotic analysis of $\hat{\pi}_{\nu,N,M}^{\mathrm{SIR}}(f)$ is quite delicate because of the dependencies in the SIR sample (which, again, is not an i.i.d. sample from $\pi$).

2.3. *The population Monte Carlo algorithm.* In their generalization of importance sampling, Cappé et al. [3] introduce an iterative dimension in the production of importance samples, aimed at adapting the importance distribution $\nu$ to the target distribution $\pi$. Iterations are then used to learn about $\pi$ from the (poor or good) performance of earlier proposals and that performance can be evaluated using different criteria, as, for example, the entropy of the distribution of importance weights.

More precisely, at iteration $t$ $(t = 0, 1, \ldots, T)$ of the PMC algorithm, a sample of $N$ values from the target distribution $\pi$ is produced by a SIR algorithm whose importance function $\nu_t$ depends on $t$, in the sense that $\nu_t$ can be derived from the $N \times (t-1)$ previous realizations of the algorithm, except for the first iteration $t = 0$, where the proposal distribution $\nu_0$ is chosen as in a regular importance sampling experiment. A generic rendering of the algorithm is thus as follows:

PMC ALGORITHM.  At time $t = 0, 1, \ldots, T$,

1. generate $(x_{i,t})_{1 \leq i \leq N}$ by i.i.d. sampling from $\nu_t$ and compute the normalized importance weights $\bar{\omega}_{i,t}$;
2. resample $(\tilde{x}_{i,t})_{1 \leq i \leq N}$ from $(x_{i,t})_{1 \leq i \leq N}$ by multinomial sampling with weights $\bar{\omega}_{i,t}$;
3. construct $\nu_{t+1}$ based on $(x_{i,\tau}, \bar{\omega}_{i,\tau})_{1 \leq i \leq N, 0 \leq \tau \leq t}$.

At this stage of the description of the PMC algorithm, we do not give in detail the construction of $\nu_t$, which can thus arbitrarily depend on the past simulations. Section 3 studies a particular choice of the $\nu_t$'s in detail. A major finding of Cappé et al. [3] is, however, that the dependence of $\nu_t$ on earlier proposals and realizations does not jeopardize the fundamental importance sampling identity. Local adaptive importance sampling schemes can thus be chosen in much wider generality than was previously thought and this shows that, thanks to the introduction of a temporal dimension in the selection of the importance function, an adaptive perspective can be adopted at little cost, with a potentially large gain in efficiency.

Obviously, when the construction of the proposal distribution $\nu_t$ is completely open, there is no guarantee of permanent improvement of the simulation scheme, whatever the criterion adopted to measure this improvement.

For instance, as illustrated by Cappé et al. [3], a constant dependence of $\nu_t$ on the past sample quickly leads to stable behavior. A more extreme illustration is the case where the sequence $(\nu_t)$ degenerates into a quasi-Dirac mass at a value based on earlier simulations and where the performance of the resulting PMC algorithm worsens. In order to study the positive and negative effects of the update of the importance function $\nu_t$, we hereafter consider a special class of proposals based on mixtures and study two particular updating schemes, one in which improvement does not occur and one in which it does.

**3. The $D$-kernel PMC algorithm.** In this section and the following ones, we study adaptivity for a particular type of parameterized PMC scheme, in the case where $\nu_t$ is a mixture of measures $\zeta_d$ ($1 \leq d \leq D$) that are chosen prior to the simulation experiment, based either on an educated guess about $\pi$ or on local approximations (as in nonparametric kernel estimation). The dependence of the $\nu_t$'s on the past simulations is of the form

$$\nu_t(dx) = \sum_{d=1}^{D} \alpha_d^t \zeta_d(\{x_{i,t-1}, \bar{\omega}_{i,t-1}\}_{1 \leq i \leq N}, dx)$$

$$= \sum_{d=1}^{D} \alpha_d^t \sum_{j=1}^{N} \bar{\omega}_{j,t-1} Q_d(x_{j,t-1}, dx),$$

where the $\bar{\omega}_{i,t}$'s denote the importance weights and the transition kernels $Q_d$ ($1 \leq d \leq D$) are given. This situation is rather common in MCMC settings where several competing transition kernels are simultaneously available, but difficult to compare. For instance, the cycle and mixture MCMC schemes discussed by Tierney [25] are of this nature.

Hereafter, we thus focus on adapting the weights $(a_d^t)_{1 \leq d \leq D}$ toward a better fit with the target distribution. A natural approach to updating the weights $(a_d^t)_{1 \leq d \leq D}$ is to favor those kernels $Q_d$ that lead to a high acceptance probability in the resampling step of the PMC algorithm. We thus choose the $a_d^t$'s to be proportional to the survival rates of the corresponding $Q_d$'s in the previous resampling step [since using the double mixture of the measures $Q_d(x_{j,t-1}, dx)$ means that we first select a point $x_{j,t-1}$ in the previous sample with probability $\bar{\omega}_{i,t-1}$, then select a component $d$ with probability $\alpha_d^t$ and, finally, simulate from $Q_d(x_{j,t-1}, dx)$].

3.1. *Algorithmic details.* The family $(Q_d)_{1 \leq d \leq D}$ of transition kernels on $\Omega \times \mathcal{A}$ is such that $(Q_d(x, \cdot))_{1 \leq d \leq D, \, x \in \Omega}$ is dominated by the reference measure $\mu$ introduced earlier. As above, we set the corresponding target density and the transition kernel to be $\pi$ and $q_d(\cdot, \cdot)$ ($1 \leq d \leq D$), respectively.

The associated PMC algorithm then updates the proposal weights (and generates the corresponding samples from $\pi$) as follows:

$D$-KERNEL PMC ALGORITHM.    At time 0,

(a) generate $x_{i,0} \overset{\text{i.i.d.}}{\sim} \nu_0$ $(1 \leq i \leq N)$ and compute the normalized importance weights $\bar{\omega}_{i,0} \propto \{\pi/\nu_0\}(x_{i,0})$;

(b) resample $(x_{i,0})_{1 \leq i \leq N}$ into $(\tilde{x}_{i,0})_{1 \leq i \leq N}$ by multinomial sampling

$$(J_{i,0})_{1 \leq i \leq N} \sim \mathcal{M}(N, (\bar{\omega}_{i,0})_{1 \leq i \leq N}), \quad \text{that is,} \quad \tilde{x}_{i,0} = x_{J_{i,0},0};$$

(c) set $\alpha_d^{1,N} = 1/D$ $(1 \leq d \leq D)$.

At time $t = 1, \ldots, T$,

(a) select the mixture components $(K_{i,t})_{1 \leq i \leq N} \sim \mathcal{M}(N, (\alpha_d^{t,N})_{1 \leq d \leq D})$;

(b) generate independent $x_{i,t} \sim q_{K_{i,t}}(\tilde{x}_{i,t-1}, x)$ $(1 \leq i \leq N)$ and compute the normalized importance weights $\bar{\omega}_{i,t} \propto \pi(x_{i,t})/q_{K_{i,t}}(\tilde{x}_{i,t-1}, x_{i,t})$;

(c) resample $(x_{i,t})_{1 \leq i \leq N}$ into $(\tilde{x}_{i,t})_{1 \leq i \leq N}$ by multinomial sampling with weights $(\bar{\omega}_{i,t})_{1 \leq i \leq N}$;

(d) set $\alpha_d^{t+1,N} = \sum_{i=1}^{N} \bar{\omega}_{i,t} \mathbb{I}_d(K_{i,t})$ $(1 \leq d \leq D)$.

In this implementation, at time $t \geq 1$, step (a) chooses a kernel index $d$ in the mixture for each point of the sample, while step (d) updates the weight $\alpha_d$ as the relative importance of kernel $Q_d$ in the current round or, in other words, as the relative *survival rate* of the points simulated from kernel $Q_d$. (Indeed, since the survival of a simulated value $x_{i,t}$ is driven by its importance weight $\bar{\omega}_{i,t}$, reweighting is related to the respective magnitudes of the importance weights for the different kernels.) Also, note that the resampling step (c) is used to avoid the propagation of very small importance weights along iterations and the subsequent degeneracy phenomenon that plagues iterated IS schemes like particle filters [11]. At time $t = T$, the algorithm should thus stop at step (b) for integral approximations to be based on the weighted sample $(x_{i,T}, \bar{\omega}_{i,T})_{1 \leq i \leq N}$.

3.2. *Convergence properties.*    In order to assess the impact of this update mechanism on the performance of the PMC scheme, we now consider the convergence properties of the above algorithm when the sample size $N$ grows to infinity. Indeed, as already pointed out in Cappé et al. [3], it does not make sense to consider the alternative asymptotics of the PMC scheme (namely, when $T$ grows to infinity), given that this algorithm is intended to be run with a small number $T$ of iterations.

A basic assumption on the kernels $Q_d$ is that the corresponding importance weights are almost surely finite, that is,

(A1) $\qquad \forall d \in \{1, \ldots, D\} \qquad \pi \otimes \pi \{q_d(x, x') = 0\} = 0,$

where $\xi \otimes \zeta$ denotes the product measure, that is, $\xi \otimes \zeta(A \times B) = \int_{A \times B} \xi(dx) \times \zeta(dy)$. We denote by $\gamma_u$ the uniform distribution on $\{1, \ldots, D\}$, that is, $\gamma_u(k) = 1/D$ for all $k \in \{1, \ldots, D\}$. The following result (whose proof is given in [Appendix B]) is a general LLN on the pairs $(x_{i,t}, K_{i,t})$ produced by the above algorithm.

PROPOSITION 3.1. *Under* [(A1)], *for any* $\pi \otimes \gamma_u$-*measurable function* $h$ *and every* $t \in \mathbb{N}$,

$$\sum_{i=1}^{N} \bar{\omega}_{i,t} h(x_{i,t}, K_{i,t}) \xrightarrow[\mathbb{P}]{N \to \infty} \pi \otimes \gamma_u(h).$$

Note that this convergence result is more than we need for Monte Carlo purposes since the $K_{i,t}$'s are auxiliary variables that are not relevant to the original problem. However, the consequences of this general result are far from negligible. First, it implies that the approximation provided by the algorithm is convergent, in the sense that $\sum_{i=1}^{N} \bar{\omega}_{i,t} f(x_{i,t})$ is a convergent estimator of $\pi(f)$. But, more importantly, it also shows that for $t \geq 1$,

$$\alpha_d^{t,N} = \sum_{i=1}^{N} \bar{\omega}_{i,t} \mathbb{I}_d(K_{i,t}) \xrightarrow[\mathbb{P}]{N \to \infty} 1/D.$$

Therefore, at *each* iteration, the weights of *all* kernels converge to $1/D$ when the sample size grows to infinity. This translates into a lack of learning properties for the $D$-kernel PMC algorithm: its properties at iteration 1 and at iteration 10 are the same. In other words, this algorithm is not adaptive and only requires one iteration for a large value of $N$. We can also relate this to the fast stabilization of the approximation in [3]. (Note that a CLT can be established for this algorithm, but given its unappealing features, we leave the exercise for the interested reader.)

**4. The Kullback divergence.** In order to obtain a correct and adaptive version of the $D$-kernel algorithm, we must first choose an effective criterion to evaluate both the adaptivity and the approximation of the target distribution by the proposal distribution. We then propose a modification of the original $D$-kernel algorithm that achieves efficiency in this sense.

As argued in many papers, using a wide range of arguments, a natural choice of approximation metric is the Kullback divergence. We thus aim to derive the $D$-kernel mixture that minimizes the Kullback divergence between this mixture and the target measure $\pi$,

$$(1) \qquad \iint \log\left(\frac{\pi(x)\pi(x')}{\pi(x)\sum_{d=1}^{D} \alpha_d q_d(x, x')}\right)(\pi \otimes \pi)(dx, dx').$$

This section is devoted to the problem of finding an iterative choice of mixing coefficients $\alpha_d$ that converge to this minimum. A detailed description of the optimal PMC scheme is given in Section [5].

4.1. *The criterion.* Using the same notation as above, in conjunction with the choice of weights $\alpha_d$ in the $D$-kernel mixture, we introduce the simplex of $\mathbb{R}^D$,

$$\mathscr{S} = \left\{ \alpha = (\alpha_1, \ldots, \alpha_D); \ \forall d \in \{1, \ldots, D\}, \alpha_d \geq 0 \text{ and } \sum_{d=1}^{D} \alpha_d = 1 \right\},$$

and denote $\bar{\pi} = \pi \otimes \pi$. We now assume that the $D$ kernels also satisfy the condition

(A2)
$$\forall d \in \{1, \ldots, D\}$$
$$\mathbb{E}_{\bar{\pi}}[|\log q_d(X, X')|] = \iint |\log q_d(x, x')| \bar{\pi}(dx, dx') < \infty,$$

which is automatically satisfied when all ratios $\pi/q_d$ are bounded. From the Kullback divergence, we derive a function on $\mathscr{S}$ such that for $\alpha \in \mathscr{S}$,

$$\mathcal{E}_{\bar{\pi}}(\alpha) = \iint \log \sum_{d=1}^{D} \alpha_d q_d(x, x') \bar{\pi}(dx, dx') = \mathbb{E}_{\bar{\pi}} \left[ \log \sum_{d=1}^{D} \alpha_d q_d(X, X') \right].$$

By virtue of Jensen's inequality, $\mathcal{E}_{\bar{\pi}}(\alpha) \leq \int \pi(dx) \log \pi(x)$ for all $\alpha \in \mathscr{S}$.

Note that due to the strict concavity of the log function, $\mathcal{E}_{\bar{\pi}}$ is a strictly concave function on a connected compact set and thus has no local maximum besides the global maximum, denoted $\alpha^{\max}$. Since

$$\int \log \pi(x) \pi(dx) - \mathcal{E}_{\bar{\pi}}(\alpha) = \mathbb{E}_{\bar{\pi}} \left( \log \pi(X) \pi(X') \middle/ \pi(X) \left\{ \sum_{d=1}^{D} \alpha_d q_d(X, X') \right\} \right),$$

$\alpha^{\max}$ is the optimal vector of weights for a mixture of transition kernels such that the product of $\pi$ by this mixture is the closest to the product distribution $\bar{\pi}$.

4.2. *A maximization algorithm.* We now study an iterative procedure, akin to the EM algorithm, that updates the weights so that the function $\mathcal{E}_{\bar{\pi}}(\alpha)$ increases at each step. Defining the function $F$ on $\mathscr{S}$ as

$$F(\alpha) = \left( \mathbb{E}_{\bar{\pi}} \left[ \alpha_d q_d(X, X') \middle/ \sum_{j=1}^{D} \alpha_j q_j(X, X') \right] \right)_{1 \leq d \leq D},$$

we construct the sequence $(\alpha^t)_{t \geq 1}$ on $\mathscr{S}$ such that $\alpha^1 = (1/D, \ldots, 1/D)$ and $\alpha^{t+1} = F(\alpha^t)$ for $t \geq 1$. Note that under assumption (A1), for all $t \geq 0$,

$$\mathbb{E}_{\bar{\pi}} \left( q_d(X, X') \middle/ \sum_{j=1}^{D} \alpha_j^t q_j(X, X') \right) > 0$$

and thus for all $t \geq 0$ and $1 \leq d \leq D$, we have $\alpha_d^t > 0$. If we define the extremal set

$$\mathcal{I}_D = \left\{ \alpha \in \mathscr{S}; \forall d \in \{1, \ldots, D\}, \text{ either } \alpha_d = 0 \text{ or } \right.$$

$$\left. \mathbb{E}_{\bar{\pi}} \left( \frac{q_d(X, X')}{\sum_{j=1}^{D} \alpha_j q_j(X, X')} \right) = 1 \right\},$$

we then have the following fixed-point result:

PROPOSITION 4.1. *Under* (A1) *and* (A2),

(i) $\mathcal{E}_{\bar{\pi}} \circ F - \mathcal{E}_{\bar{\pi}}$ *is continuous;*
(ii) *for all* $\alpha \in \mathscr{S}$, $\mathcal{E}_{\bar{\pi}} \circ F(\alpha) \geq \mathcal{E}_{\bar{\pi}}(\alpha)$;
(iii) $\mathcal{I}_D = \{\alpha \in \mathscr{S}; F(\alpha) = \alpha\} = \{\alpha \in \mathscr{S}; \mathcal{E}_{\bar{\pi}} \circ F(\alpha) = \mathcal{E}_{\bar{\pi}}(\alpha)\}$ *and* $\mathcal{I}_D$ *is finite.*

PROOF. $\mathcal{E}_{\bar{\pi}}$ is clearly continuous. Moreover, by Lebesgue's dominated convergence theorem, the function $\alpha \mapsto \mathbb{E}_{\bar{\pi}}(\alpha_d q_d(X, X')/\sum_{j=1}^{d} \alpha_j q_j(X, X'))$ is also continuous, which implies that $F$ is continuous. This completes the proof of (i). Due to the concavity of the log function,

$$\mathcal{E}_{\bar{\pi}}(F(\alpha)) - \mathcal{E}_{\bar{\pi}}(\alpha)$$

$$= \mathbb{E}_{\bar{\pi}} \left( \log \left[ \sum_{d=1}^{D} \frac{\alpha_d q_d(X, X')}{\sum_{j=1}^{D} \alpha_j q_j(X, X')} \mathbb{E}_{\bar{\pi}} \left( \frac{q_d(X, X')}{\sum_{j=1}^{D} \alpha_j q_j(X, X')} \right) \right] \right)$$

$$\geq \mathbb{E}_{\bar{\pi}} \left[ \sum_{d=1}^{D} \frac{\alpha_d q_d(X, X')}{\sum_{j=1}^{D} \alpha_j q_j(X, X')} \log \mathbb{E}_{\bar{\pi}} \left( \frac{q_d(X, X')}{\sum_{j=1}^{D} \alpha_j q_j(X, X')} \right) \right]$$

$$= \sum_{d=1}^{D} \alpha_d \mathbb{E}_{\bar{\pi}} \left( \frac{q_d(X, X')}{\sum_{j=1}^{D} \alpha_j q_j(X, X')} \right) \log \mathbb{E}_{\bar{\pi}} \left( \frac{q_d(X, X')}{\sum_{j=1}^{D} \alpha_j q_j(X, X')} \right).$$

Applying the inequality $u \log u \geq u - 1$ yields (ii). Moreover, the equality in $u \log u \geq u - 1$ holds if and only if $u = 1$. Therefore, equality above is equivalent to

$$\forall \alpha_d \neq 0 \qquad \mathbb{E}_{\bar{\pi}} \left( q_d(X, X') \Big/ \sum_{j=1}^{D} \alpha_j q_j(X, X') \right) = 1.$$

Thus, $\mathcal{I}_D = \{\alpha \in \mathscr{S}; \mathcal{E}_{\bar{\pi}} \circ F(\alpha) = \mathcal{E}_{\bar{\pi}}(\alpha)\}$. The second equality, $\mathcal{I}_D = \{\alpha \in \mathscr{S}; F(\alpha) = \alpha\}$, is straightforward.

We now prove by recursion on $D$ that $\mathcal{I}_D$ is finite, which is equivalent to proving that

$$\{\alpha \in \mathcal{I}_D; \ \alpha_d \neq 0 \ \forall d \in \{1, \ldots, D\}\}$$

is empty or finite. If this set is nonempty, then any element $\alpha$ in this set satisfies

$$\forall d \in \{1, \ldots, D\} \qquad \mathbb{E}_{\bar{\pi}}\left(\frac{q_d(X, X')}{\sum_{j=1}^{D} \alpha_j q_j(X, X')}\right) = 1,$$

which implies that

$$\begin{aligned}
0 &= \sum_{d=1}^{D} \alpha_d^{\max}\left(\mathbb{E}_{\bar{\pi}}\left(\frac{q_d(X, X')}{\sum_{j=1}^{D} \alpha_j q_j(X, X')}\right) - 1\right) \\
&= \mathbb{E}_{\bar{\pi}}\left(\frac{\sum_{d=1}^{D} \alpha_d^{\max} q_d(X, X')}{\sum_{j=1}^{D} \alpha_j q_j(X, X')} - 1\right) \\
&\geq \mathbb{E}_{\bar{\pi}}\left(\log \frac{\sum_{d=1}^{D} \alpha_d^{\max} q_d(X, X')}{\sum_{j=1}^{D} \alpha_j q_j(X, X')}\right) \geq 0.
\end{aligned}$$

Since the global maximum of $\mathcal{E}_{\bar{\pi}}$ is unique, we conclude that $\alpha = \alpha^{\max}$ and hence (iii) follows. $\square$

4.3. *Averaged EM.* Proposition 4.1 implies that our recursive procedure satisfies $\mathcal{E}_{\bar{\pi}}(\alpha^{t+1}) \geq \mathcal{E}_{\bar{\pi}}(\alpha^t)$. Therefore, the Kullback divergence criterion (1) decreases at each step. This property is closely linked to the EM algorithm ([20], Section 5.3). More precisely, consider the mixture model

$$V \sim \mathcal{M}(1, \{\alpha_1, \ldots, \alpha_D\}) \quad \text{and} \quad W = (X, X') | V \sim \pi(dx) Q_V(x, dx')$$

with parameter $\alpha$. We denote by $\bar{\bar{\mathbb{E}}}_\alpha$ the corresponding expectation, by $p_\alpha(v, w)$ the joint density of $(V, W)$ with respect to $\mu \otimes \mu$ and by $p_\alpha(w)$ the density of $W$ with respect to $\mu$. It is then easy to check that $\mathcal{E}_{\bar{\pi}}(\alpha) = \int \log(p_\alpha(w)) \bar{\pi}(dw)$, which is an average version of the criterion to be maximized in the EM algorithm when only $W$ is observed. In this case, a natural idea, adapted from the EM algorithm, is to update $\alpha$ according to the iterative scheme

$$\alpha^{t+1} = \arg\max_{\alpha \in \mathscr{S}} \int \bar{\bar{\mathbb{E}}}_{\alpha^t}[\log p_\alpha(V, w) | w] \bar{\pi}(dw).$$

Straightforward algebra can be used to show that this definition of $\alpha^{t+1}$ is equivalent to the update formula $\alpha^{t+1} = F(\alpha^t)$ that we used above. Our algorithm then appears as an averaged EM, but shares with EM the property that the criterion increases deterministically at each step.

The following result ensures that any $\alpha$ different from $\alpha^{\max}$ is repulsive.

PROPOSITION 4.2. *Under* (A1) *and* (A2), *for every* $\alpha \neq \alpha^{\max} \in \mathscr{S}$, *there exists a neighborhood* $V_\alpha$ *of* $\alpha$ *such that if* $\alpha^{t_0} \in V_\alpha$, *then* $(\alpha^t)_{t \geq t_0}$ *leaves* $V_\alpha$ *within a finite time.*

PROOF.    Let $\alpha \neq \alpha^{\max}$. Then using the inequality $u - 1 \geq \log u$, we have

$$\sum_{d=1}^{D} \alpha_d^{\max} \mathbb{E}_{\bar{\pi}} \left( \frac{q_d(X, X')}{\sum_{j=1}^{D} \alpha_j q_j(X, X')} \right) - 1$$

$$\geq \mathbb{E}_{\bar{\pi}} \left( \log \frac{\sum_{d=1}^{D} \alpha_d^{\max} q_d(X, X')}{\sum_{j=1}^{D} \alpha_j q_j(X, X')} \right) > 0,$$

which implies that there exists $1 \leq d \leq D$ such that

$$\mathbb{E}_{\bar{\pi}} \left( q_d(X, X') \Big/ \sum_{j=1}^{D} \alpha_j q_j(X, X') \right) > 1.$$

Using a nonincreasing sequence $(W_n)_{n \geq 0}$ of neighborhoods of $\alpha$ in $\mathscr{S}$, the monotone convergence theorem implies that

$$1 < \mathbb{E}_{\bar{\pi}} \left( \frac{q_d(X, X')}{\sum_{j=1}^{D} \alpha_j q_j(X, X')} \right) = \mathbb{E}_{\bar{\pi}} \left( \lim_{n \to \infty} \inf_{\beta \in W_n} \frac{q_d(X, X')}{\sum_{j=1}^{D} \beta_j q_j(X, X')} \right)$$

$$\leq \lim_{n \to \infty} \inf_{\beta \in W_n} \mathbb{E}_{\bar{\pi}} \left( \frac{q_d(X, X')}{\sum_{j=1}^{D} \beta_j q_j(X, X')} \right).$$

Thus, there exist $W_{n_0} = V_\alpha$, a neighborhood of $\alpha$ and $\eta > 1$ such that

$$(2) \qquad \forall \beta \in V_\alpha \qquad \mathbb{E}_{\bar{\pi}} \left( q_d(X, X') \Big/ \sum_{j=1}^{D} \beta_j q_j(X, X') \right) > \eta.$$

Now, for all $t \geq 0$ and $1 \leq d \leq D$, we have $1 \geq \alpha_d^t > 0$. Combining (2) with the update formulas for $\alpha^t = F(\alpha^{t-1})$ shows that $(\alpha^t)_{t \geq 0}$ leaves $V_\alpha$ within a finite time.    □

We thus conclude that the maximization algorithm is convergent, as asserted by the following proposition:

PROPOSITION 4.3.    *Under* (A1) *and* (A2),

$$\lim_{t \to \infty} \alpha^t = \alpha^{\max}.$$

PROOF.    First, recall that $\mathcal{I}_D$ is a finite set and that $\alpha^{\max} \in \mathcal{I}_D$, that is, $\mathcal{I}_D = \{\beta_0, \beta_1, \dots, \beta_I\}$ with $\beta_0 = \alpha^{\max}$. We introduce a sequence $(W_i)_{0 \leq i \leq I}$ of disjoint neighborhoods of the $\beta_i$'s such that for all $0 \leq i \leq I$, $F(W_i)$ is disjoint from $\bigcup_{j \neq i} W_j$ [this is possible since $F(\beta_i) = \beta_i$ and $F$ is continuous] and for all $i \in \{1, \dots, I\}$, $W_i \subset V_{\beta_i}$, where the $(V_{\beta_i})$'s are defined as in the proof of Proposition 4.2.

The sequence $(\mathcal{E}_{\bar{\pi}}(\alpha^t))_{t \geq 0}$ is upper-bounded and nondecreasing; therefore it converges. This implies that $\lim_{t \to \infty} \mathcal{E}_{\bar{\pi}} \circ F(\alpha^t) - \mathcal{E}_{\bar{\pi}}(\alpha^t) = 0$. By continuity of $\mathcal{E}_{\bar{\pi}} \circ F - \mathcal{E}_{\bar{\pi}}$, there exists $T > 0$ such that for all $t \geq T$, $\alpha_t \in \bigcup_j W_j$. Since $F(W_i)$ is disjoint from $\bigcup_{j \neq i} W_j$, this implies that there exists $i \in \{0, \dots, I\}$ such that for all $t \geq T$, $\alpha^t \in W_i$. By Proposition 4.2, $i$ cannot be in $\{1, \dots, I\}$. Thus, for all $t \geq T$, $\alpha^t \in W_0$, which is a neighborhood of $\beta_0 = \alpha^{\max}$.   $\square$

## 5. The Rao–Blackwellized $D$-kernel PMC.

Update of the weights $\alpha_d$ through the transform $F$ thus improves the Kullback divergence criterion. We now discuss how to implement this mechanism within a PMC algorithm that resembles the previous $D$-kernel algorithm. The only difference with the algorithm of Section 3.1 is that we make use of the kernel structure in the computation of the importance weight. In MCMC terminology, this is called "Rao–Blackwellization" ([20], Section 4.2) and it is known to provide variance reduction in data augmentation settings ([20], Section 9.2). In the current context, the improvement brought about by Rao–Blackwellization is dramatic, in that the modified algorithm does converge to the proposal mixture that is closest to the target distribution in the sense of the Kullback divergence. More precisely, a Monte Carlo version of the update via $F$ can be implemented in the iterative definition of the mixture weights, in the same way that MCEM approximates EM ([20], Section 5.3.3).

5.1. *The algorithm.* In importance sampling, as well as in MCMC settings, the (de)conditioning improvement brought about by Rao–Blackwellization may be significant [5]. In the case of the $D$-kernel PMC scheme, the Rao–Blackwellization argument is that it is not necessary to condition on the value of the mixture component in the computation of the importance weight and that the improvement is brought about by using the whole mixture. The importance weight should therefore be

$$\pi(x_{i,t}) \bigg/ \sum_{d=1}^{D} \alpha_d^{t,N} q_d(\tilde{x}_{i,t-1}, x_{i,t}) \quad \text{rather than} \quad \pi(x_{i,t})/q_{K_{i,t}}(\tilde{x}_{i,t-1}, x_{i,t}),$$

as in the algorithm of Section 3.1. As already noted by Hesterberg [15], the use of the whole mixture in the importance weight is a robust tool that prevents infinite variance importance sampling estimators. In the next section, we show that this choice of weight guarantees that the following modification of the $D$-kernel algorithm converges to the optimal mixture (in terms of Kullback divergence).

RAO–BLACKWELLIZED $D$-KERNEL PMC ALGORITHM. At time 0, use the same steps as in the $D$-kernel PMC algorithm to obtain $(\tilde{x}_{i,0})_{1 \leq i \leq N}$ and set $\alpha_d^{1,N} = 1/D$ $(1 \leq d \leq D)$.

At time $t = 1, \dots, T$,

(a) generate $(K_{i,t})_{1 \leq i \leq N} \sim \mathcal{M}(N, (\alpha_d^{t,N})_{1 \leq d \leq D})$;

(b) generate independent $x_{i,t} \sim q_{K_{i,t}}(\tilde{x}_{i,t-1}, x)$ $(1 \le i \le N)$ and compute the normalized importance weights $\bar{\omega}_{i,t} \propto \pi(x_{i,t})/\sum_{d=1}^{D} \alpha_d^{t,N} q_d(\tilde{x}_{i,t-1}, x_{i,t})$;

(c) resample $(x_{i,t})$ into $(\tilde{x}_{i,t})_{1 \le i \le N}$ by multinomial sampling with weights $(\bar{\omega}_{i,t})_{1 \le i \le N}$;

(d) set $\alpha_d^{t+1,N} = \sum_{i=1}^{N} \bar{\omega}_{i,t} \mathbb{I}_d(K_{i,t})$ $(1 \le d \le D)$.

5.2. *The corresponding law of large numbers.* Not very surprisingly, the sample obtained at each iteration of the above Rao–Blackwellized algorithm approximates the target distribution in the sense of the weak law of large numbers (LLN). Note that the convergence holds under the very weak assumption (A1) and for any test function $h$ that is absolutely integrable with respect to the target distribution $\pi$. The function $h$ may thus be unbounded.

THEOREM 5.1. *Under* (A1), *for any function $h$ in $L_\pi^1$ and for all $t \ge 0$,*

$$\sum_{i=1}^{N} \bar{\omega}_{i,t} h(x_{i,t}) \overset{N \to \infty}{\underset{\mathbb{P}}{\longrightarrow}} \pi(h) \quad and \quad \frac{1}{N} \sum_{i=1}^{N} h(\tilde{x}_{i,t}) \overset{N \to \infty}{\underset{\mathbb{P}}{\longrightarrow}} \pi(h).$$

PROOF. First, convergence of the second average follows from the convergence of the first average by Theorem A.1, since the latter is simply a multinomial sampling perturbation of the former. We thus focus on the weighted average and proceed by induction on $t$. The case $t = 0$ is the basic importance sampling convergence result. For $t \ge 1$, if the convergence holds for $t - 1$, then to prove convergence at iteration $t$, we need only check, as in Proposition 3.1, that

$$\frac{1}{N} \sum_{i=1}^{N} \omega_{i,t} h(x_{i,t}) \overset{N \to \infty}{\underset{\mathbb{P}}{\longrightarrow}} \pi(h),$$

where $\omega_{i,t}$ denotes the importance weight $\pi(x_{i,t})/\sum_{d=1}^{D} \alpha_d^{t,N} q_d(\tilde{x}_{i,t-1}, x_{i,t})$. (The special case $h \equiv 1$ ensures that the renormalizing sum converges to 1.) We apply Theorem A.1 with $\mathcal{G}_N = \sigma((\tilde{x}_{i,t-1})_{1 \le i \le N}, (\alpha_d^{t,N})_{1 \le d \le D})$ and $U_{N,i} = N^{-1} \omega_{i,t} h(x_{i,t})$. Then conditionally on $\mathcal{G}_N$, the $x_{i,t}$'s $(1 \le i \le N)$ are independent and

$$x_{i,t} | \mathcal{G}_N \sim \sum_{d=1}^{D} \alpha_d^{t,N} Q_d(\tilde{x}_{i,t-1}, \cdot).$$

Noting that

$$\sum_{i=1}^{N} \mathbb{E}\left( \frac{\omega_{i,t} h(x_{i,t})}{N} \bigg| \mathcal{G}_N \right)$$

$$= \sum_{i=1}^{N} \mathbb{E}\left( \frac{\pi(x_{i,t}) h(x_{i,t})}{N \sum_{d=1}^{D} \alpha_d^{t,N} q_d(\tilde{x}_{i,t-1}, x_{i,t})} \bigg| \mathcal{G}_N \right) = \pi(h),$$

we only need to check condition (iii). The end of the proof is then quite similar to the proof of Proposition 3.1.  □

5.3. *Convergence of the weights.* The next proposition ensures that at each iteration, the update of the mixture weights in the Rao–Blackwellized algorithm approximates the theoretical update obtained in Section 4.2 for minimizing the Kullback divergence criterion.

PROPOSITION 5.1.  *Under* (A1), *for all* $t \geq 1$,

$$\forall 1 \leq d \leq D \qquad \alpha_d^{t,N} \xrightarrow[\mathbb{P}]{N \to \infty} \alpha_d^t, \tag{3}$$

*where* $\alpha^t = F(\alpha^{t-1})$.

Combining Proposition 5.1 with Proposition 4.3, we obtain that, under assumptions (A1)–(A2), the Rao–Blackwellized version of the PMC algorithm adapts the weights of the proposed mixture of kernels, in the sense that it converges to the optimal combination of mixtures with respect to the Kullback divergence criterion obtained in Section 4.2.

PROOF OF PROPOSITION 5.1.  The case $t = 1$ is obvious. Now, assume (3) holds for some $t \geq 1$. As in the proof of Proposition 3.1, we now establish that

$$\frac{1}{N} \sum_{i=1}^N \omega_{i,t} \mathbb{I}_d(K_{i,t}) = \frac{1}{N} \sum_{i=1}^N \frac{\pi(x_{i,t})}{\sum_{l=1}^D \alpha_l^{t,N} q_l(\tilde{x}_{i,t-1}, x_{i,t})} \mathbb{I}_d(K_{i,t}) \xrightarrow[\mathbb{P}]{N \to \infty} \alpha_d^{t+1},$$

the convergence of the renormalizing sum to 1 being a consequence of this convergence. We apply Theorem A.1 with $\mathcal{G}_N = \sigma((\tilde{x}_{i,t-1})_{1 \leq i \leq N}, (\alpha_d^{t,N})_{1 \leq d \leq D})$ and $U_{N,i} = N^{-1} \omega_{i,t} \mathbb{I}_d(K_{i,t})$. Conditionally on $\mathcal{G}_N$, the $(K_{i,t}, x_{i,t})$'s $(1 \leq i \leq N)$ are independent and for all $(d, A)$ in $\{1, \ldots, D\} \times \mathcal{A}$, we have

$$\mathbb{P}(K_{i,t} = d, x_{i,t} \in A | \mathcal{G}_N) = \alpha_d^{t,N} Q_d(\tilde{x}_{i,t-1}, A).$$

To apply Theorem A.1, we need only check condition (iii). We have, for $C > 0$,

$$\mathbb{E}\left( \sum_{i=1}^N \frac{\omega_{i,t} \mathbb{I}_d(K_{i,t})}{N} \mathbb{I}_{\{\omega_{i,t} \mathbb{I}_d(K_{i,t}) > C\}} \,\bigg|\, \mathcal{G}_N \right)$$

$$\leq \sum_{j=1}^D \frac{1}{N} \sum_{i=1}^N \int \frac{\alpha_d^{t,N} q_d(\tilde{x}_{i,t-1}, x)}{\sum_{l=1}^D \alpha_l^{t,N} q_l(\tilde{x}_{i,t-1}, x)} \mathbb{I}\left\{ \frac{\pi(x)}{D^{-1} q_j(\tilde{x}_{i,t-1}, x)} > C \right\} \pi(dx)$$

$$\leq \sum_{j=1}^D \frac{1}{N} \sum_{i=1}^N \int \mathbb{I}\left\{ \frac{\pi(x)}{D^{-1} q_j(\tilde{x}_{i,t-1}, x)} > C \right\} \pi(dx)$$

$$\xrightarrow[\mathbb{P}]{N\to\infty} \sum_{j=1}^{D} \bar{\pi}\left(\frac{\pi(x)}{D^{-1}q_j(x',x)} > C\right),$$

by the LLN of Theorem 5.1. The right-hand side converges to 0 as $C$ tends to infinity since by assumption (A1), $\bar{\pi}(\{q_j(x,x') = 0\}) = 0$. Thus, Theorem A.1 applies and

$$\frac{1}{N}\sum_{i=1}^{N}\omega_{i,t}\mathbb{I}_d(K_{i,t}) - \mathbb{E}\left(\frac{1}{N}\sum_{i=1}^{N}\omega_{i,t}\mathbb{I}_d(K_{i,t})\Big|\mathcal{G}_N\right) \xrightarrow[\mathbb{P}]{N\to\infty} 0.$$

To complete the proof, it simply remains to show that

$$
\begin{aligned}
(4) \qquad &\mathbb{E}\left(\frac{1}{N}\sum_{i=1}^{N}\omega_{i,t}\mathbb{I}_d(K_{i,t})\Big|\mathcal{G}_N\right)\\
&= \frac{1}{N}\sum_{i=1}^{N}\int \frac{\alpha_d^{t,N}q_d(\tilde{x}_{i,t-1},x)}{\sum_{l=1}^{D}\alpha_l^{t,N}q_l(\tilde{x}_{i,t-1},x)}\pi(dx) \xrightarrow[\mathbb{P}]{N\to\infty} \alpha_d^{t+1}.
\end{aligned}
$$

It follows from the LLN stated in Theorem 5.1 that

$$(5) \qquad \frac{1}{N}\sum_{i=1}^{N}\int \pi(dx)\frac{\alpha_d^t q_d(\tilde{x}_{i,t-1},x)}{\sum_{l=1}^{D}\alpha_l^t q_l(\tilde{x}_{i,t-1},x)} \xrightarrow[\mathbb{P}]{N\to\infty} \mathbb{E}_{\bar{\pi}}\left(\frac{\alpha_d^t q_d(X,X')}{\sum_{l=1}^{D}\alpha_l^t q_l(X,X')}\right) = \alpha_d^{t+1}$$

and it thus suffices to check that the difference between (4) and (5) converges to 0 in probability. To show this, first note that for all $t \geq 1$ and all $1 \leq d \leq D$, $\alpha_d^t > 0$ and thus, by the induction assumption, for all $1 \leq d \leq D$, $(\alpha_d^{t,N} - \alpha_d^t)/\alpha_d^t \xrightarrow[\mathbb{P}]{N\to\infty} 0$. Using the inequality $|\frac{A}{B} - \frac{C}{D}| \leq |\frac{A}{B}|\frac{D-B}{D}| + |\frac{A-C}{C}||\frac{C}{D}|$, we have, by straightforward algebra, that

$$
\begin{aligned}
&\left|\frac{\alpha_d^{t,N}q_d(\tilde{x}_{i,t-1},x)}{\sum_{l=1}^{D}\alpha_l^{t,N}q_l(\tilde{x}_{i,t-1},x)} - \frac{\alpha_d^t q_d(\tilde{x}_{i,t-1},x)}{\sum_{l=1}^{D}\alpha_l^t q_l(\tilde{x}_{i,t-1},x)}\right|\\
&\qquad \leq \frac{\alpha_d^{t,N}q_d(\tilde{x}_{i,t-1},x)}{\sum_{j=1}^{D}\alpha_l^{t,N}q_j(\tilde{x}_{i,t-1},x)}\left(\sup_{l\in\{1,\ldots,D\}}\left|\frac{\alpha_l^{t,N} - \alpha_l^t}{\alpha_l^t}\right|\right)\\
&\qquad\quad + \left|\frac{\alpha_d^{t,N} - \alpha_d^t}{\alpha_d^t}\right|\left|\frac{\alpha_d^t q_d(\tilde{x}_{i,t-1},x)}{\sum_{l=1}^{D}\alpha_l^t q_l(\tilde{x}_{i,t-1},x)}\right|\\
&\qquad \leq 2\sup_{l\in\{1,\ldots,D\}}\left|\frac{\alpha_l^{t,N} - \alpha_l^t}{\alpha_l^t}\right|.
\end{aligned}
$$

The proposition then follows from the convergence $(\alpha_d^{t,N} - \alpha_d^t)/\alpha_d^t \xrightarrow[\mathbb{P}]{N\to\infty} 0$. $\square$

5.4. *A corresponding central limit theorem.* We now establish a CLT for the weighted and the unweighted samples when the sample size goes to infinity. As noted in Section 2.2 for the SIR algorithm, the asymptotic variance associated with the unweighted sample is larger than the variance of the weighted sample because of the additional multinomial step.

THEOREM 5.2. *Under* (A1),

(i) *for a function h such that* $\bar{\pi}\{h^2(x')\pi(x)/q_d(x,x')\} < \infty$ *for at least one* $1 \le d \le D$, *we have*

$$(6) \qquad \sqrt{N} \sum_{i=1}^{N} \bar{\omega}_{i,t}\{h(x_{i,t}) - \pi(h)\} \xrightarrow{\mathcal{L}} \mathcal{N}(0, \sigma_t^2),$$

*where* $\sigma_t^2 = \bar{\pi}(\{h(x') - \pi(h)\}^2 \pi(x')/\sum_{d=1}^{D} \alpha_d^t q_d(x,x'))$;

(ii) *if, moreover,* $\pi(h^2) < \infty$, *then*

$$(7) \qquad \frac{1}{\sqrt{N}} \sum_{i=1}^{N}\{h(\tilde{x}_{i,t}) - \pi(h)\} \xrightarrow{\mathcal{L}} \mathcal{N}(0, \sigma_t^2 + \mathbb{V}_{\pi}(h)).$$

Note that amongst the conditions under which this theorem applies, the integrability condition

$$(8) \qquad \bar{\pi}\left(h^2(x')\frac{\pi(x)}{q_d(x,x')}\right) < \infty$$

is required for *some d* in $\{1, \dots, D\}$ and not for *all d*'s. Thus, settings where some transition kernels $q_d(\cdot, \cdot)$ do not satisfy (8) can still be covered by this theorem provided (8) holds for at least one particular kernel. An equivalent expression of the asymptotic variance $\sigma_t^2$ is

$$\sigma_t^2 = \mathbb{V}_{\nu}\left(\{h - \pi(h)\}\frac{\bar{\pi}}{\nu}\right) \qquad \text{where } \nu(dx, dx') = \pi(dx)\left(\sum_{d=1}^{D} \alpha_d^t Q_d(x, dx')\right).$$

Written in this form, $\sigma_t^2$ turns out to be the expression of the asymptotic variance that appears in the CLT associated with a self-normalized IS algorithm (SNIS) (see Section 2.1 for a description of the algorithm), where the proposal distribution would be $\nu$ and the target distribution $\bar{\pi}$. Obviously, this SNIS algorithm cannot be implemented since, given the above definition of $\nu$, the proposal distribution depends on both $\pi$ and the weights $(\alpha_d^t)$, which are unknown.

PROOF OF THEOREM 5.2. Without loss of generality, we assume that $\pi(h) = 0$. Let $d_0 \in \{1, \dots, D\}$ be such that $\bar{\pi}(h^2(x')\pi(x)/q_{d_0}(x,x')) < \infty$.

A consequence of the proof of Theorem 5.1 is that $\frac{1}{N}\sum_{i=1}^{N}\omega_{i,t} \xrightarrow[\mathbb{P}]{N\to\infty} 1$, so we only need to prove that

$$(9) \qquad \frac{1}{\sqrt{N}}\sum_{i=1}^{N}\omega_{i,t}h(x_{i,t}) \xrightarrow{\mathcal{L}} \mathcal{N}(0,\sigma_t^2).$$

We will apply Theorem A.2 with

$$U_{N,i} = \frac{1}{\sqrt{N}}\omega_{i,t}h(x_{i,t}) = \frac{1}{\sqrt{N}}\frac{\pi(x_{i,t})h(x_{i,t})}{\sum_{d=1}^{D}\alpha_d^{t,N}q_d(\tilde{x}_{i,t-1},x_{i,t})},$$

$$\mathcal{G}_N = \sigma\{(\tilde{x}_{i,t-1})_{1\le i\le N}, (\alpha_d^{t,N})_{1\le d\le D}\}.$$

Conditionally on $\mathcal{G}_N$, the $(x_{i,t})$'s $(1\le i\le N)$ are independent and

$$x_{i,t}|\mathcal{G}_N \sim \sum_{d=1}^{D}\alpha_d^{t,N}Q_d(\tilde{x}_{i,t-1},\cdot).$$

Conditions (i) and (ii) of Theorem A.2 are clearly satisfied. To check condition (iii), first note that $\mathbb{E}(U_{N,i}|\mathcal{G}_N) = \pi(h) = 0$. Moreover,

$$A_N = \sum_{i=1}^{N}\mathbb{E}(U_{N,i}^2|\mathcal{G}_N)$$

$$= \frac{1}{N}\sum_{i=1}^{N}\int h^2(x)\frac{\pi(x)}{\sum_{d=1}^{D}\alpha_d^{t,N}q_d(\tilde{x}_{i,t-1},x)}\pi(dx).$$

By the LLN for $(\tilde{x}_{i,t})$ stated in Theorem 5.1, we have

$$B_N = \frac{1}{N}\sum_{i=1}^{N}\int h^2(x)\frac{\pi(x)}{\sum_{d=1}^{D}\alpha_d^{t}q_d(\tilde{x}_{i,t-1},x)}\pi(dx) \xrightarrow[\mathbb{P}]{N\to\infty} \sigma_t^2.$$

To prove that (iii) holds, it is thus sufficient to show that $|B_N - A_N| \xrightarrow[\mathbb{P}]{N\to\infty} 0$. Since $\alpha_{d_0}^{t,N} \xrightarrow[\mathbb{P}]{N\to\infty} \alpha_{d_0}^{t} > 0$, we need only consider the upper bound

$$\mathbb{1}_{\{\alpha_{d_0}^{t,N}>2^{-1}\alpha_{d_0}^{t}\}}|B_N - A_N|$$

$$\le \mathbb{1}_{\{\alpha_{d_0}^{t,N}>2^{-1}\alpha_{d_0}^{t}\}}\sup_{1\le d\le D}\left(\frac{\alpha_d^{t}-\alpha_d^{t,N}}{\alpha_d^{t}}\right)\frac{1}{N}\sum_{j=1}^{N}\int \frac{h^2(x)\pi(x)}{\sum_{d=1}^{D}\alpha_d^{t,N}q_d(\tilde{x}_{i,t-1},x)}\pi(dx)$$

$$\le \left\{\sup_{1\le d\le D}\left(\frac{\alpha_d^{t}-\alpha_d^{t,N}}{\alpha_d^{t}}\right)\frac{1}{N}\sum_{j=1}^{N}\int \frac{h^2(x)\pi(x)}{2^{-1}\alpha_{d_0}^{t}q_{d_0}(\tilde{x}_{i,t-1},x)}\right\}\pi(dx) \xrightarrow[\mathbb{P}]{N\to\infty} 0.$$

Thus, condition (iii) is satisfied. Finally, we consider condition (iv). Using the same argument as was used for condition (iii), we have that

$$
\mathbb{I}_{\{\alpha_{d_0}^{t,N} > 2^{-1}\alpha_{d_0}^t\}} \sum_{i=1}^{N} \mathbb{E}\left[\frac{1}{N}\omega_{i,t}^2 h^2(x_{i,t})\mathbb{I}\left\{\frac{\pi(x_{i,t})h(x_{i,t})}{\sum_{d=1}^{D}\alpha_d^{t,N}q_d(\tilde{x}_{i,t-1},x_{i,t})} > C\right\}\Big|\mathcal{G}_N\right]
$$

$$
\leq \frac{1}{N}\sum_{i=1}^{N}\int h^2(x)\frac{\pi(x)}{2^{-1}\alpha_{d_0}^t q_{d_0}(\tilde{x}_{i,t-1},x)}\mathbb{I}\left\{\frac{\pi(x)h(x)}{2^{-1}\alpha_{d_0}^t q_{d_0}(\tilde{x}_{i,t-1},x)} > C\right\}\pi(dx)
$$

$$
\xrightarrow[N\to\infty]{\mathbb{P}}\bar{\pi}\left(h^2(x)\frac{\pi(x)}{2^{-1}\alpha_{d_0}^t q_{d_0}(x',x)}\mathbb{I}\left\{\frac{\pi(x)h(x)}{2^{-1}\alpha_{d_0}^t q_{d_0}(x',x)} > C\right\}\right),
$$

which converges to 0 as $C$ tends to infinity. Thus, Theorem A.2 can be applied and the proof of (6) is completed. The proof of (7) follows from a direct application of Theorem A.2, as in the SIR result, by setting $U_{N,i} = \frac{1}{\sqrt{N}}h(\tilde{x}_{i,t})$ and $\mathcal{G}_N = \sigma\{(x_{i,t})_{1\leq i\leq N},(\omega_{i,t})_{1\leq i\leq N}\}$.  □

**6. Illustrations.** In this section, we briefly show how the iterations of the PMC algorithm quickly implement adaptivity toward the most efficient mixture of kernels, through three examples of moderate difficulty. (The R programs are available from the last author's website via an Snw file.)

Example 1. As a first toy example, we take the target $\pi$ to be the density of a five-dimensional normal mixture,

$$
(10) \qquad\qquad \sum_{i=1}^{3}\frac{1}{3}\mathcal{N}_5(0,\Sigma_i),
$$

and an independent normal mixture proposal with the same means and variances as (10), but started with different weights $\alpha_d^{2,N}$. Note that this is a very special case of a $D$-kernel PMC scheme in that the transition kernels of Section 5.1 are then independent proposals. In this case, the optimal choice of weights is obviously $\alpha_d^\star = 1/3$. In our experiment, we used three Wishart-simulated variances $\Sigma_i$, with 10, 15 and 7 degrees of freedom, respectively. The starting values $\alpha_d^{1,N}$ are indicated on the left of Figure 1, which clearly shows the convergence to the optimal values 1/3 and 2/3 for the two first accumulated weights in less than ten iterations. (Generating more simulated points at each iteration stabilizes the convergence graph, but the primary aim of this example is to exhibit the fast convergence to the true optimal values of the weights.)

EXAMPLE 2. As a second toy example, consider the case of a three-dimensional normal $\mathscr{N}_3(0, \Sigma)$ target with covariance matrix

$$\Sigma = \begin{pmatrix} 6.986 & 0.154 & 3.523 \\ 0.154 & 15.433 & 3.528 \\ 3.523 & 3.528 & 18.463 \end{pmatrix}$$

and the mixture of three kernels given by

$$(11) \qquad \alpha_1^{t,N} \mathscr{T}_2(\tilde{x}_{i,t-1}, \Sigma_1) + \alpha_2^{t,N} \mathscr{N}(\tilde{x}_{i,t-1}, \Sigma_2) + \alpha_3^{t,N} \mathscr{N}(\tilde{x}_{i,t-1}, \Sigma_3),$$

where the $\Sigma_i$'s are random Wishart matrices. [The first proposal in the mixture is a product of unidimensional Student $t$-distributions with two degrees of freedom, centered at the current values of the components of $\tilde{x}_{i,t-1}$ and rotated by $\sqrt{\tau_1} \operatorname{diag}(\Sigma^{1/2})$.]

A compelling feature of this example is that we can visualize the Kullback divergence on the $\mathbb{R}^3$ simplex in the sense that the divergence

$$\mathfrak{E}(\pi, \tilde{\pi}) = \mathbb{E}_{\tilde{\pi}} \left[ \log \left( \pi(X') \middle/ \left\{ \sum_{d=1}^{3} \alpha_d q_d(X, X') \right\} \right) \right] = \mathfrak{e}(\alpha_1, \alpha_2, \alpha_3)$$

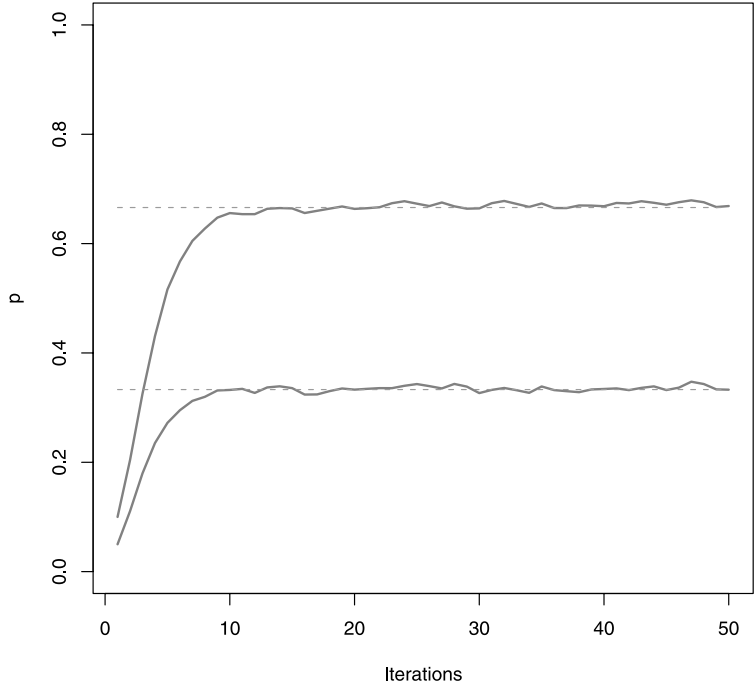

FIG. 1. *Convergence of the accumulated weights $\alpha_1^{t,N}$ and $\alpha_1^{t,N} + \alpha_2^{t,N}$ for the three-component normal mixture to the optimal values 1/3 and 2/3 (represented by dotted lines). At each iteration, $N = 1{,}000$ points were simulated from the D-kernel proposal.*

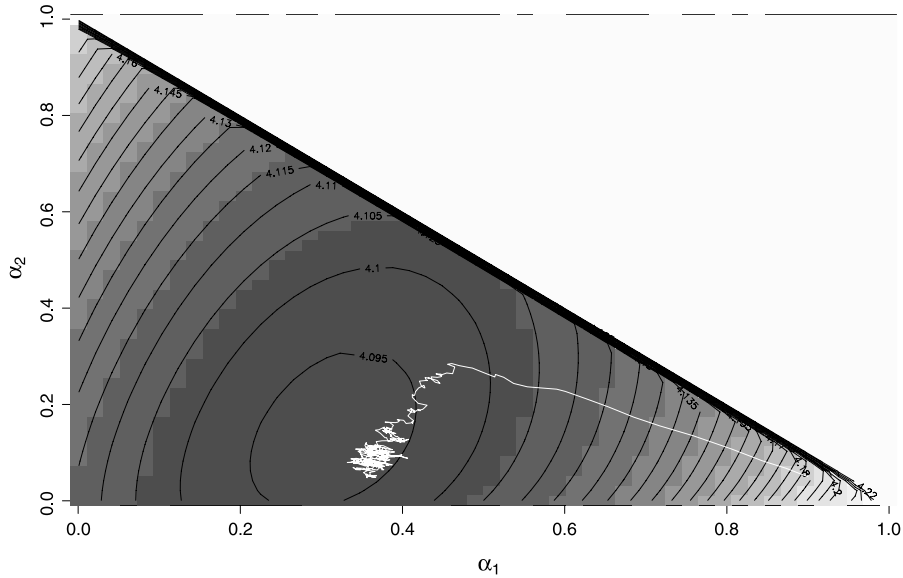

FIG. 2.   *Grey level and contour representation of the Kullback divergence $\mathfrak{k}(\alpha_1, \alpha_2, \alpha_3)$ between the $\mathscr{N}_3(0, \Sigma)$ distribution and the three-component mixture proposal* (11). (*The darker pixels correspond to lower values of the divergence.*) *We also represent* (*in white*) *the path of one run of the D-kernel PMC algorithm when started from a random value* $(\alpha_1, \alpha_2, \alpha_3)$. *The number of iterations T is equal to* 500, *while the sample size N is* 50,000.

can be approximated by a Monte Carlo experiment on a grid of values of $(\alpha_1, \alpha_2)$. Figure 2 shows the result of this Monte Carlo experiment based on 25,000 $\mathscr{N}_3(0, \Sigma)$ simulations and exhibits a minimum divergence inside the $\mathbb{R}^3$ simplex. Running the Rao–Blackwellized $D$-kernel PMC algorithm from a random starting weight $(\alpha_1, \alpha_2, \alpha_3)$ always leads to a neighborhood of the minimum, even though a strict decrease in the divergence requires a large value for $N$ and a precise convergence to $\alpha^{\max}$ necessitates a large number of iterations $T$.

EXAMPLE 3.   Our third example is a contingency table inspired by [1], given here as Table 1. We model this dataset by a Poisson regression,

$$x_{ij} \sim \mathscr{P}(\exp(\alpha_i + \beta_j)), \qquad i, j = 0, 1,$$

with $\alpha_0 = 0$ for identifiability reasons. We use a flat prior on the parameter $\theta = (\alpha_1, \beta_0, \beta_1)$ and run the PMC $D$-kernel algorithm with a mixture of ten normal random walk proposals, $\mathscr{N}(\tilde{\theta}_{i,t-1}, \varrho_d I(\hat{\theta}))$, $d = 1, \ldots, 10$, where $I(\hat{\theta})$ is the Fisher information matrix evaluated at the MLE, $\hat{\theta} = (-0.43, 4.06, 5.9)$, and where the scales $\varrho_d$ vary from $1.35e{-}19$ to $1.54e{+}07$ (the $\varrho_d$'s are equidistributed on a logarithmic scale). The results of five (successive) iterations of the Rao–Blackwell $D$-kernel algorithm are as follows: unsurprisingly, the largest variance kernels are hardly ever sampled, but fulfill their main role of variance stabilizers in the importance sampling weights while the mixture concentrates on the medium variances, with a quick convergence of the mixture weights to the limiting weights—the accumulated weights of the 5th, 6th, 7th and 8th components of the mixture converge to 0, 0.003, 0.259 and 0.738, respectively. The fit of the simulated sample to the target distribution is shown in Figure 3, since the points of the sample do coincide with the (unique) modal region of the posterior distribution. This experiment also shows that there is no degeneracy in the samples produced: most points in the last sample have very similar posterior values. For instance, 20% of the sample corresponds to 95% of the weights, while 1% of the sample corresponds to 31% of the weights. A closer look at convergence is provided by Figure 4, where the histograms of the resampled samples are represented, along with the distribution of the log likelihood and the empirical cumulative distribution function cdf of the importance weights. They do not signal any degeneracy phenomenon, but, rather the opposite—a clear stabilization around the values of interest.

**7. Conclusion and perspectives.** This paper shows that it is possible to build an adaptive mixture of proposals aimed at a minimization of the Kullback divergence with the distribution of interest. We can therefore set different goals for a simulation experiment and expect to arrive at the most accurate proposal. For instance, in a companion paper [9], we also derive an adaptive update of the weights targeted at the minimal variance proposal for a given integral $\mathfrak{I}$ of interest. Rather naturally, these results are achieved under strong restrictions on the family of proposals in the sense that the parameterization of those families is restricted to the weights of the mixture. It is, however, possible to extend the above results to general mixtures of parameterized proposals, as shown by work currently under development.

TABLE 1
*Two-by-two contingency table*

|  | **0** | **1** | **Total** |
|---|---|---|---|
| 0 | 60 | 364 | 424 |
| 1 | 36 | 240 | 276 |
| Total | 96 | 604 | 700 |

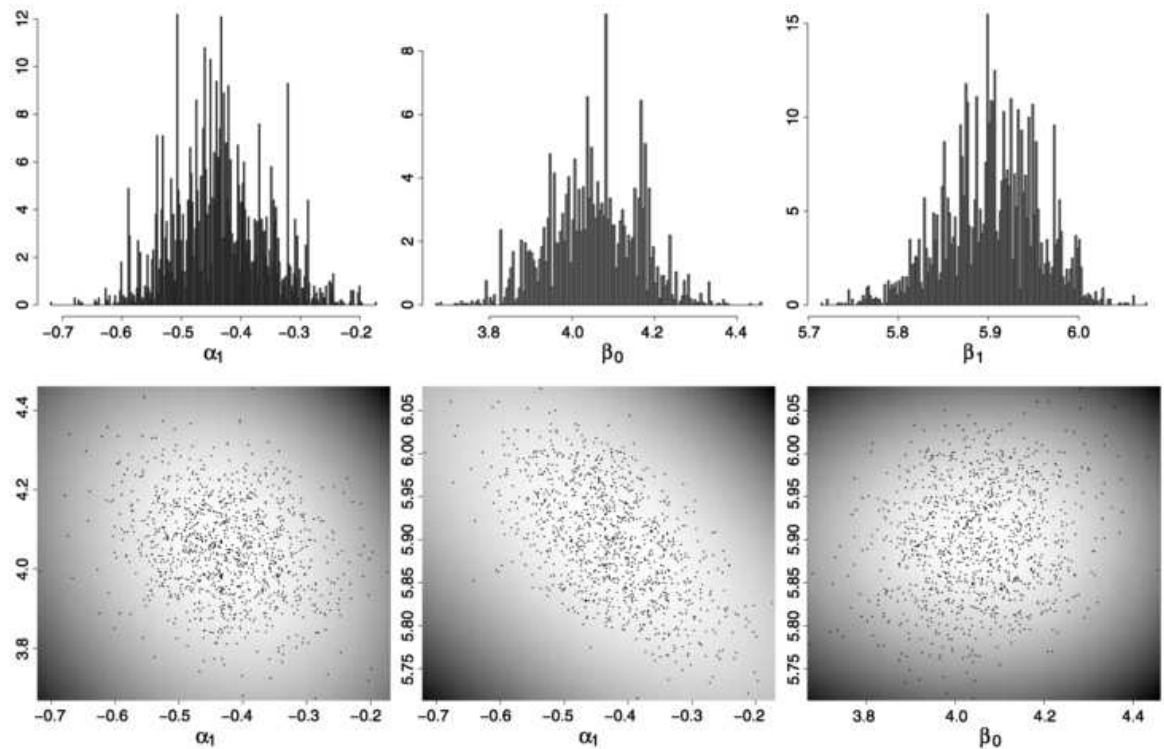

FIG. 3.  *Distribution of* 5,000 *resampled points after five iterations of the Rao–Blackwellized D-kernel PMC sampler for the contingency table example.* Top: *histograms of the components* $\alpha_1$, $\beta_0$ *and* $\beta_1$; bottom: *scatterplots of the points* $(\alpha_1, \beta_0)$, $(\alpha_1, \beta_1)$ *and* $(\beta_0, \beta_1)$ *on the profile slices of the log-likelihood.*

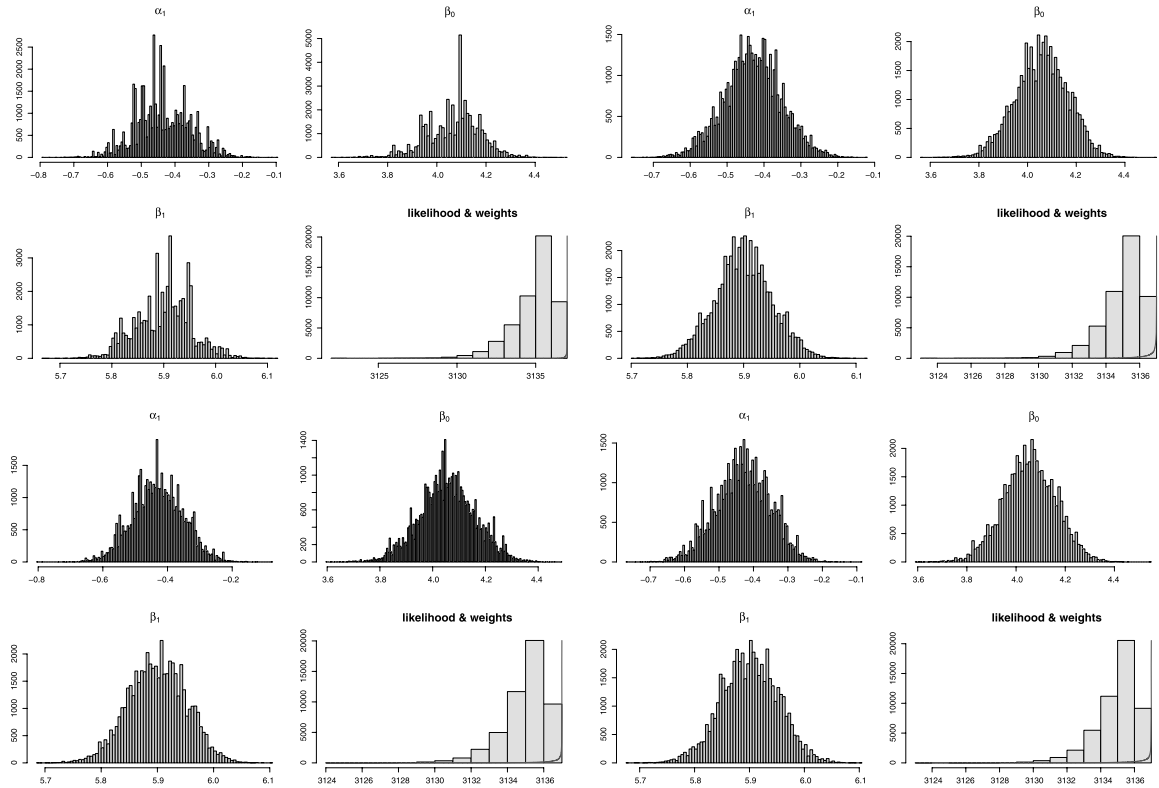

FIG. 4. *Evolution of the samples over four iterations of the Rao–Blackwellized D-kernel PMC sampler for the contingency table example (the output from each iteration is a block of four graphs, to be read from left to right and from top to bottom): histograms of the resampled samples of $\alpha_1$, $\beta_0$ and $\beta_1$ of size 50,000 and (lower right of each block) log-likelihood and the empirical cumulative distribution function of the importance weights.*

A more practical direction of research is the implementation of such adaptive algorithms in large-dimensional problems. While our algorithms are in fine importance sampling algorithms, it is conceivable that mixtures of Gibbs-like proposals can take better advantage of the intuition gained from MCMC methodology, while keeping the finite-horizon validation of importance sampling methods. The major difficulty in this direction, however, is that the curse of dimensionality still holds in the sense that (a) we need to simultaneously consider more and more proposals as the dimension increases (as, e.g., the set of all full conditionals) and (b) the number of parameters to tune in the proposals exponentially increases with the dimension.

## APPENDIX A: CONVERGENCE THEOREMS FOR TRIANGULAR ARRAYS

In this section, we recall convergence results for triangular arrays of random variables (see [4] or [10] for more details, including the proofs). We will use these results to study the asymptotic behavior of the PMC algorithm. In what follows, let $\{U_{N,i}\}_{N \geq 1, 1 \leq i \leq N}$ be a triangular array of random variables defined on the same measurable space $(\Omega, \mathcal{A})$ and let $\{\mathcal{G}_N\}_{N \geq 1}$ be a sequence of $\sigma$-algebras included in $\mathcal{A}$. The symbol $X_N \longrightarrow_{\mathbb{P}} a$ means that $X_N$ converges in probability to $a$ as $N$ goes to infinity.

DEFINITION A.1.   The sequence $\{U_{N,i}\}_{N \geq 1, 1 \leq i \leq N}$ is *independent given* $\{\mathcal{G}_N\}_{N \geq 1}$ if for all $N \geq 1$, the random variables $U_{N,1}, \ldots, U_{N,N}$ are independent given $\mathcal{G}_N$.

DEFINITION A.2.   The sequence of variables $\{Z_N\}_{N \geq 1}$ is *bounded in probability* if

$$\lim_{C \to \infty} \sup_{N \geq 1} \mathbb{P}[|Z_N| \geq C] = 0.$$

THEOREM A.1.   *If*

(i)  $\{U_{N,i}\}_{N \geq 1, 1 \leq i \leq N}$ *is independent given* $\{\mathcal{G}_N\}_{N \geq 1}$;
(ii)  *the sequence* $\{\sum_{i=1}^{N} \mathbb{E}[|U_{N,i}||\mathcal{G}_N]\}_{N \geq 1}$ *is bounded in probability*;
(iii)  *for all* $\eta > 0$, $\sum_{i=1}^{N} \mathbb{E}[|U_{N,i}|\mathbb{I}_{|U_{N,i}| > \eta}|\mathcal{G}_N] \longrightarrow_{\mathbb{P}} 0$,

*then* $\sum_{i=1}^{N}(U_{N,i} - \mathbb{E}[U_{N,i}|\mathcal{G}_N]) \longrightarrow_{\mathbb{P}} 0$.

THEOREM A.2.   *If*

(i)  $\{U_{N,i}\}_{N \geq 1, 1 \leq i \leq N}$ *is independent given* $\{\mathcal{G}_N\}_{N \geq 1}$;
(ii)  *for all* $N \geq 1, \forall i \in \{1, \ldots, N\}$, $\mathbb{E}[|U_{N,i}||\mathcal{G}_N] < \infty$;

(iii) *there exists* $\sigma^2 > 0$ *such that* $\sum_{i=1}^{N}(\mathbb{E}[U_{N,i}^2|\mathcal{G}_N] - (\mathbb{E}[U_{N,i}|\mathcal{G}_N])^2) \longrightarrow_{\mathbb{P}} \sigma^2$;

(iv) *for all* $\eta > 0$, $\sum_{i=1}^{N}\mathbb{E}[U_{N,i}^2\mathbb{I}_{|U_{N,i}|>\eta}|\mathcal{G}_N] \longrightarrow_{\mathbb{P}} 0$,

*then for all* $u \in \mathbb{R}$,

$$\mathbb{E}\left[\exp\left(iu\sum_{i=1}^{N}(U_{N,i} - \mathbb{E}[U_{N,i}|\mathcal{G}_N])\right)\Big|\mathcal{G}_N\right] \longrightarrow_{\mathbb{P}} \exp\left(-\frac{u^2\sigma^2}{2}\right).$$

## APPENDIX B: PROOF OF PROPOSITION 3.1

We proceed by induction with respect to $t$. Using Theorem A.1, the case $t = 0$ is straightforward as a direct consequence of the convergence of the importance sampling algorithm.

For $t \geq 1$, let us assume that the LLN holds for $t-1$. Then to prove that $\sum_{i=1}^{N}\bar{\omega}_{i,t}h(x_{i,t}, K_{i,t})$ converges in probability to $\pi \otimes \gamma_u(h)$, we need only check that

$$N^{-1}\sum_{i=1}^{N}\frac{\pi(x_{i,t})}{q_{K_{i,t}}(\tilde{x}_{i,t-1}, x_{i,t})}h(x_{i,t}, K_{i,t}) \xrightarrow[\mathbb{P}]{N\to\infty} \pi \otimes \gamma_u(h),$$

the special case $h \equiv 1$ providing the convergence of the normalizing constant for the importance weights. Applying Theorem A.1 with

$$U_{N,i} = N^{-1}\frac{\pi(x_{i,t})}{q_{K_{i,t}}(\tilde{x}_{i,t-1}, x_{i,t})}h(x_{i,t}, K_{i,t})$$

and

$$\mathcal{G}_N = \sigma\{(\tilde{x}_{i,t-1})_{1\leq i\leq N}, (\alpha_d^{t,N})_{1\leq d\leq D}\},$$

where $\sigma\{(X_i)_i\}$ denotes the $\sigma$-algebra induced by the $X_i$'s, we need only check condition (iii). For any $C > 0$, we have

$$
\begin{aligned}
N^{-1}\sum_{i=1}^{N}\mathbb{E}\bigg[&\frac{\pi(x_{i,t})}{q_{K_{i,t}}(\tilde{x}_{i,t-1}, x_{i,t})} \\
&\times h(x_{i,t}, K_{i,t})\mathbb{I}\bigg\{\frac{\pi(x_{i,t})}{q_{K_{i,t}}(\tilde{x}_{i,t-1}, x_{i,t})}h(x_{i,t}, K_{i,t}) > C\bigg\}\bigg|\mathcal{G}_N\bigg] \\
&= \sum_{d=1}^{D}N^{-1}\sum_{i=1}^{N}F_C(\tilde{x}_{i,t-1}, d)\alpha_d^{t,N},
\end{aligned}
$$
(12)

where $F_C(x, k) = \int \pi(du)h(u, k)\mathbb{I}\{\frac{\pi(u)}{q_k(x,u)}h(u, k) \geq C\}$. By induction, we have

$$\alpha_d^{t,N} = \sum_{i=1}^{N} \bar{\omega}_{i,t-1}\mathbb{I}_d(K_{i,t-1}) \underset{\mathbb{P}}{\longrightarrow} 1/D \quad \text{and}$$

$$N^{-1} \sum_{i=1}^{N} F_C(\tilde{x}_{i,t-1}, k) \overset{N\to\infty}{\underset{\mathbb{P}}{\longrightarrow}} \pi(F_C(\cdot, k)).$$

Using these limits in (12) yields

$$N^{-1} \sum_{i=1}^{N} \mathbb{E}\left[\frac{\pi(x_{i,t})h(x_{i,t}, K_{i,t})}{q_{K_{i,t}}(\tilde{x}_{i,t-1}, x_{i,t})}\mathbb{I}\left\{\frac{\pi(x_{i,t})h(x_{i,t}, K_{i,t})}{q_{K_{i,t}}(\tilde{x}_{i,t-1}, x_{i,t})} > C\right\}\bigg|\mathcal{G}_N\right]$$

$$\overset{N\to\infty}{\underset{\mathbb{P}}{\longrightarrow}} \pi \otimes \gamma_u(F_C).$$

Since $\pi \otimes \gamma_u(F_C)$ converges to 0 as $C$ goes to infinity, this proves that for any $\eta > 0$,

$$N^{-1} \sum_{i=1}^{N} \mathbb{E}\left[\frac{\pi(x_{i,t})h(x_{i,t}, K_{i,t})}{q_{K_{i,t}}(\tilde{x}_{i,t-1}, x_{i,t})}\mathbb{I}\left\{\frac{\pi(x_{i,t})h(x_{i,t}, K_{i,t})}{q_{K_{i,t}}(\tilde{x}_{i,t-1}, x_{i,t})} > N\eta\right\}\bigg|\mathcal{G}_N\right] \overset{N\to\infty}{\underset{\mathbb{P}}{\longrightarrow}} 0.$$

Condition (iii) is satisfied and Theorem A.1 applies. The proof follows. □

**Acknowledgments.** The authors are grateful to Olivier Cappé, Paul Fearnhead and Eric Moulines for helpful comments and discussions. Comments from two referees helped considerably in improving the focus and presentation of our results.

## REFERENCES

[1] AGRESTI, A. (2002). *Categorical Data Analysis*, 2nd ed. Wiley, New York. MR1914507

[2] ANDRIEU, C. and ROBERT, C. (2001). Controlled Markov chain Monte Carlo methods for optimal sampling. Technical Report 0125, Univ. Paris Dauphine.

[3] CAPPÉ, O., GUILLIN, A., MARIN, J. and ROBERT, C. (2004). Population Monte Carlo. *J. Comput. Graph. Statist.* **13** 907–929. MR2109057

[4] CAPPÉ, O., MOULINES, E. and RYDÉN, T. (2005). *Inference in Hidden Markov Models.* Springer, New York. MR2159833

[5] CELEUX, G., MARIN, J. and ROBERT, C. (2006). Iterated importance sampling in missing data problems. *Comput. Statist. Data Anal.* **50** 3386–3404.

[6] CHOPIN, N. (2004). Central limit theorem for sequential Monte Carlo methods and its application to Bayesian inference. *Ann. Statist.* **32** 2385–2411. MR2153989

[7] CSISZÁR, I. and TUSNÁDY, G. (1984). Information geometry and alternating minimization procedures. Recent results in estimation theory and related topics. *Statist. Decisions* **1984** (suppl. 1) 205–237. MR0785210

[8] DEL MORAL, P., DOUCET, A. and JASRA, A. (2006). Sequential Monte Carlo samplers. *J. R. Stat. Soc. Ser. B Stat. Methodol.* **68** 411–436. MR2278333

[9] DOUC, R., GUILLIN, A., MARIN, J. and ROBERT, C. (2005). Minimum variance importance sampling via population Monte Carlo. Technical report, Cahiers du CEREMADE, Univ. Paris Dauphine.

[10] DOUC, R. and MOULINES, E. (2005). Limit theorems for properly weighted samples with applications to sequential Monte Carlo. Technical report, TSI, Telecom Paris.

[11] DOUCET, A., DE FREITAS, N. and GORDON, N., eds. (2001). *Sequential Monte Carlo Methods in Practice.* Springer, New York. MR1847783

[12] GILKS, W., ROBERTS, G. and SAHU, S. (1998). Adaptive Markov chain Monte Carlo through regeneration. *J. Amer. Statist. Assoc.* **93** 1045–1054. MR1649199

[13] HAARIO, H., SAKSMAN, E. and TAMMINEN, J. (1999). Adaptive proposal distribution for random walk Metropolis algorithm. *Comput. Statist.* **14** 375–395.

[14] HAARIO, H., SAKSMAN, E. and TAMMINEN, J. (2001). An adaptive Metropolis algorithm. *Bernoulli* **7** 223–242. MR1828504

[15] HESTERBERG, T. (1995). Weighted average importance sampling and defensive mixture distributions. *Technometrics* **37** 185–194.

[16] IBA, Y. (2000). Population-based Monte Carlo algorithms. *Trans. Japanese Society for Artificial Intelligence* **16** 279–286.

[17] KÜNSCH, H. (2005). Recursive Monte Carlo filters: Algorithms and theoretical analysis. *Ann. Statist.* **33** 1983–2021. MR2211077

[18] MENGERSEN, K. L. and TWEEDIE, R. L. (1996). Rates of convergence of the Hastings and Metropolis algorithms. *Ann. Statist.* **24** 101–121. MR1389882

[19] ROBERT, C. (1996). Intrinsic losses. *Theory and Decision* **40** 191–214. MR1385186

[20] ROBERT, C. and CASELLA, G. (2004). *Monte Carlo Statistical Methods*, 2nd ed. Springer, New York. MR2080278

[21] ROBERTS, G. O., GELMAN, A. and GILKS, W. R. (1997). Weak convergence and optimal scaling of random walk Metropolis algorithms. *Ann. Appl. Probab.* **7** 110–120. MR1428751

[22] RUBIN, D. (1988). Using the SIR algorithm to simulate posterior distributions. In *Bayesian Statistics 3* (J. M. Bernardo, M. H. DeGroot, D. V. Lindley and A. F. M. Smith, eds.) 395–402. Oxford Univ. Press.

[23] SAHU, S. and ZHIGLJAVSKY, A. (1998). Adaptation for self regenerative MCMC. Technical report, Univ. of Wales, Cardiff.

[24] SAHU, S. and ZHIGLJAVSKY, A. (2003). Self regenerative Markov chain Monte Carlo with adaptation. *Bernoulli* **9** 395–422. MR1997490

[25] TIERNEY, L. (1994). Markov chains for exploring posterior distributions (with discussion). *Ann. Statist.* **22** 1701–1762. MR1329166

R. DOUC
CMAP, ÉCOLE POLYTECHNIQUE, CNRS
ROUTE DE SACLAY
91128 PALAISEAU CEDEX
FRANCE
E-MAIL: douc@cmapx.polytechnique.fr

A. GUILLIN
ÉCOLE CENTRALE DE MARSEILLE ET LATP, CNRS
CENTRE DE MATHÉMATIQUES ET INFORMATIQUE
TECHNOPÔLE CHÂTEAU-GOMBERT
39 RUE F. JOLIOT CURIE
13453 MARSEILLE CEDEX 13
FRANCE
E-MAIL: guillin@cmi.univ-mrs.fr

J.-M. Marin
INRIA Futurs, Projet Select
Laboratoire de Mathématiques
Université d'Orsay
91405 Orsay cedex
France
E-mail: Jean-Michel.Marin@inria.fr

C. P. Robert
Ceremade, Université Paris Dauphine
75775 Paris cedex 16
France
E-mail: xian@ceremade.dauphine.fr